\documentclass{amsproc}
\usepackage{amsmath}
\usepackage{amsthm}
\usepackage{amsfonts}
\usepackage{amssymb}

\newcommand{\lan}{\langle}
\newcommand{\ra}{\rangle}

\newcommand{\g}{\mathfrak{g}}
\newcommand{\N}{\mathbb{N}}

\newcommand{\la}{{\lambda}}

\newtheorem{theorem}{Theorem}[section]
\newtheorem{lemma}[theorem]{Lemma}

\theoremstyle{definition}
\newtheorem{definition}[theorem]{Definition}
\newtheorem{example}[theorem]{Example}

\theoremstyle{remark}
\newtheorem{remark}[theorem]{Remark}

\numberwithin{equation}{section}

\begin{document}

\title{Generating functions for ranks of pre-modular categories}

\author{Eric C. Rowell}
\address{Department of Mathematics, Indiana University,
Bloomington, IN 47405}

\email{errowell@indiana.edu}

\date{August 26, 2005.}

\begin{abstract}
We derive generating functions for the ranks of pre-modular
categories associated to quantum groups at roots of unity.
\end{abstract}

\maketitle

\section{Introduction}The \emph{rank} of a semisimple tensor category $\mathcal{O}$
 is the number of simple objects (up to isomorphism).  A problem that
 has received particular attention recently is the classification
 of finite rank categories with some extra axioms satisfied by
 the representation categories of finite groups, quantum groups or Hopf algebras.
   A
\emph{pre-modular category} \cite{Br} is a finite rank ribbon
category, \emph{i.e} a semisimple, balanced, rigid braided tensor
category.  A \emph{modular category} \cite{Tur92} is a pre-modular
category that satisfies a further non-degeneracy condition (see
\cite{BK} or \cite{TuraevBook} for detailed definitions), while more
generally a \emph{fusion category} \cite{ENO} is a finite rank,
semisimple, rigid
 tensor category.  Zhenghan
Wang has recently conjectured that there are only finitely many
modular categories of a fixed rank; a conjecture that has been
verified for ranks 1,2,3 and 4, see \cite{BRSW} for an explicit list
of all modular categories of these ranks. Similar conjectures have
been proposed for pre-modular and fusion categories, and some
progress has been made for these generalizations (see \cite{Ostrik1}
and \cite{Ostrik2}).  The problem of determining the ranks of known
pre-modular categories is motivated by this conjecture and the
relationship between these categories and low-dimensional topology
\cite{TuraevBook}, quantum computing \cite{FNSWW}, \cite{FLW} and
Hopf algebras \cite{ENO}.

The most ubiquitous examples of modular categories come from two
sources: representation categories of Hopf doubles of finite group
algebras, and sub-quotients of representation categories of
quantized universal enveloping algebras of simple Lie algebras
(henceforth \emph{quantum groups}) specialized at roots of unity. In
the finite group examples one always obtains a modular category,
whereas the quantum group categories sometimes fail the modularity
condition. The purpose of this note is to give formulas for the
ranks of the pre-modular categories constructed from quantum groups.
Most of these formulas have already appeared without formal proof in
a more primitive form in \cite{survey}.  The ranks of the modular
categories constructed from finite groups have been considered in
\cite{CGR}, where all cases up to rank 50 are determined.
 \section{Pre-modular Categories from Quantum Groups}
To each simple Lie algebra $\g$ and positive integer $\ell$ (with
$\ell$ large enough) one obtains a family of pre-modular categories
sharing the same (finite) labeling set of simple objects and tensor
product decomposition rules.  The construction (which can be found
in \cite{andersen} and \cite{TuraevBook}) is summarized as follows:
Lusztig's integral form of the Drinfeld-Jimbo quantum group $U_q\g$
is well-defined for $q^2$ a primitive $\ell$th root of unity.  The
corresponding representation category is not semisimple, but has a
well-behaved subcategory of so-called \emph{tilting modules}. The
quotient of the tilting module category by the tensor ideal of
\emph{negligible} morphisms (essentially the annihilator of a
trace-form) is a pre-modular category $\mathcal{C}(\g,q,\ell)$.
 While the structure of $\mathcal{C}(\g,q,\ell)$
depends on the choice of $q$, the rank only depends on $\g$ and
$\ell$.
 To describe the labeling sets of simple objects we need some
standard notation from Lie theory, which is found in Table
\ref{notation}.
\begin{table}\label{notation}
\caption{Notation}
\begin{tabular}{*{2}{lr}}
$\la_i$ & $i$th fundamental weight\\
$P_+$ & dominant weights \\
$h$ & Coxeter number\\
$h\check{}$ & dual Coxeter number\\
$\lan\cdot,\cdot\ra$ & symmetric bilinear form on weight lattice\\
$m$ & ratio $\frac{\lan \alpha,\alpha\ra}{\lan\beta,\beta\ra}$
$\alpha$ a long root, $\beta$ a short root\\
$\theta_0$ & longest root\\
$\theta_1$ & longest short root\\
$\rho$ & half the sum of the positive roots\\
\end{tabular}
\end{table}
We take the form $\lan\cdot,\cdot\ra$ to be normalized so that the
square length of a \emph{short} root is 2.
\begin{definition}
Referring to Table \ref{notation} for notation, the
labeling sets of isomorphism classes of simple objects are (see
\cite{kirillov}):
$$C_\ell(\g):=\begin{cases}\{\la\in P_+: \lan\la+\rho,\theta_1\ra<\ell\} &
\text{if $m\mid \ell$}\\
\{\la\in P_+: \lan\la+\rho,\theta_2\ra<\ell\} & \text{if $m\nmid
\ell$}\end{cases}$$ \end{definition} Observe that $m=1$ for the
\emph{simply-laced} Lie types $A,D$ and $E$, while $m=2$ for Lie
types $B,C$ and $F_4$, and $m=3$ for Lie type $G_2$. We define an
auxiliary label $\ell_m=0$ if $m\mid\ell$ and $\ell_m=1$ if
$m\nmid\ell$ for notational convenience.

We reduce the problem of determining the cardinalities of the
labeling sets $C_\ell(\g)$ to counting partitions of $n$ with parts
in a fixed (finite) multiset $\mathcal{S}(\g,\ell_m)$ that depends
only on the rank and Lie type of $\g$ and the divisibility of $\ell$
by $m$ (encoded in $\ell_m$). Fix a simple Lie algebra $\g$ of rank
$r$ and a positive integer $\ell$. Let $\la=\sum_i a_i\la_i$ be a
dominant weight of $\g$ written as an $\N$-linear combination of
fundamental weights $\la_i$.  To determine if $\la\in C_\ell(\g)$,
we compute:
$$\lan\la+\rho,\theta_j\ra=\lan\rho,\theta_j\ra+\sum_i^r a_i\lan\la_i,\theta_j\ra$$
where $j=1$ or $2$ depending on if $m\mid\ell$ or not.  Setting
$L_i^{(j)}=\lan\la_i,\theta_j\ra$ we see that the condition that
$\la\in C_\ell(\g)$ becomes:
$$\sum_i^k a_i L_i^{(j)}\leq\ell-\lan\rho,\theta_j\ra-1.$$  Since
$a_i,L_i^{(j)}\in\N$ we have:
\begin{lemma}
The cardinality of $C_\ell(\g)$ is the number of partitions of all
natural numbers $n$, $0\leq n\leq\ell-\lan\rho,\theta_j\ra-1$ into
parts from the size $r=\rm{rank}(\g)$ multiset
$\mathcal{S}(\g,\ell_m)=[L_i^{(j)}]_i^r$.
\end{lemma}

 So it remains only to compute the numbers $\lan\rho,\theta_j\ra$
and $L_i^{(j)}$ (with $j=0,1$) for each Lie algebra $\g$ and integer
$\ell>\lan\rho,\theta_j\ra$ and to apply standard combinatorics to
count the number of partitions into parts in
$\mathcal{S}(\g,\ell_m)$. The first task is easily accomplished with
the help of the book \cite{Bou}. Table \ref{rank} lists the results
of these computations, where $\ell_0$ is the minimal non-degenerate
value of $\ell$ (\emph{i.e.} satisfying
$\ell_0\geq\lan\rho,\theta_j\ra+1$). The combinatorial techniques
are described in the next section.
\section{Generating Functions}

\begin{table}\caption{$\mathcal{C}(\g,q,\ell)$ Data}\label{rank}
\begin{tabular}{*{3}{|c}|}
\hline
$X_r$ & $\mathcal{S}(\g,\ell_m)$ & $\ell_0$\\
\hline\hline
 $A_r$ & $[1,\ldots,1]$  & $r+1$
\\ \hline  $B_r$, $\ell$ odd& $[1,2,\ldots,2]$ & $2r+1$
\\ \hline  $B_r$, $\ell$ even& $[2,2,4,\ldots,4]$  & $4r-2$
\\ \hline  $C_r$, $\ell$ odd& $[1,2,\ldots,2]$  & $2r+1$
\\ \hline  $C_r$, $\ell$ even& $[2,\ldots,2]$  & $2r+2$
\\ \hline  $D_r$& $[1,1,1,2,\ldots,2]$  & $2r-2$
\\ \hline  $E_6$& $[1,1,2,2,2,3]$ & $12$
\\ \hline  $E_7$& $[1,2,2,2,3,3,4]$ & $18$
\\ \hline  $E_8$& $[2,2,3,3,4,4,5,6]$ & $30$
\\ \hline  $F_4$, $\ell$ even& $[2,4,4,6]$ & $18$
\\ \hline  $F_4$, $\ell$ odd& $[2,2,3,4]$ & $13$
\\ \hline  $G_2$, $3 \mid \ell$& $[3,6]$ & $12$
\\ \hline  $G_2$, $3\nmid \ell$& $[2,3]$ & $7$\\
\hline
\end{tabular}
\end{table}

Let $P_{\mathcal{T}}(n)$ denote the number of partitions of $n$ into
parts in a multiset $\mathcal{T}$, and
$P_{\mathcal{T}}[s]=\sum_{n=0}^sP_{\mathcal{T}}(n)$ the number of
partitions of all integers $0\leq n\leq s$ into parts from the
multiset $\mathcal{T}$.  Any standard reference on generating
functions (see e.g. \cite{Stan}) will provide enough details about
generating functions to prove the following:
\begin{lemma}\label{gflemma}
The number $P_{\mathcal{T}}(n)$ of partitions of $n$ into parts from
the multiset $\mathcal{T}$ has generating function:
$$\prod_{t\in\mathcal{T}}\frac{1}{1-x^t}=\sum_{n=0}^\infty P_{\mathcal{T}}(n),$$
while the number $P_{\mathcal{T}}[s]$ of partitions of all $n\in\N$
with $0\leq n\leq s$ into parts from the multiset $\mathcal{T}$ has
generating function:
$$\frac{1}{1-x}\prod_{t\in\mathcal{T}}\frac{1}{1-x^t}=\sum_{n=0}^\infty P_{\mathcal{T}}[s].$$
\end{lemma}
 Applying this lemma to the sets $\mathcal{S}(\g,\ell_m)$ given in
Table \ref{rank} we obtain:
\begin{theorem}
Define
$$F_{\g,\ell_m}(x)=\frac{1}{1-x}\prod_{k\in\mathcal{S}(\g,\ell_m)}\frac{1}{1-x^k}.$$
 Then the rank $|C_\ell(\g)|$ of the pre-modular category
$\mathcal{C}(\g,q,\ell)$ is
 the coefficient of  $$x^{\ell-\ell_0+\ell_m}$$ in the power series expansion
 of $F_{\g,\ell_m}(x)$.
\end{theorem}
\begin{proof}
It is clear from Lemma \ref{gflemma} that the coefficients of
generating function $F_{\g,\ell_m}(x)$ counts the appropriate
partitions.  The coefficient of $x$ that gives the rank for a
specific $\ell$ is shifted by the minimal non-degenerate $\ell_0$,
which corresponds to the $x^0=1$ term if $m\mid\ell$ and to the
$x^1=x$ term of $m\nmid\ell$, hence the correction by $x^{\ell_m}$.
With this normalization only the coefficients of those powers of $x$
divisible (resp. indivisible)
 by $m$ give ranks corresponding to $\ell$ divisible (resp.
 indivisible) by $m$.
\end{proof}
We illustrate the application of this theorem with some examples.
\begin{example}
Let $\g$ be of type $G_2$. \begin{enumerate}
\item[(a)] Let $\ell=27$.
Then $\ell_m=0$ and $\ell_0=12$.  So the rank of
$\mathcal{C}(\g(G_2),q,27)$ is given by the $(27-12+0)=15$th
coefficient of
$$\frac{1}{(1-x)(1-x^3)(1-x^6)}=(1+x+x^2)(1+2x^3+4x^6+6x^9+9x^{12}+12x^{15}+\cdots)$$
so $|C_{27}(\g(G_2))|=15$.
\item[(b)] Let $\ell=14$.  Then $\ell_m=1$ and $\ell_0=7$.  So
$|C_{14}(\g(G_2))|$ is the $(14-7+1)$th coefficient of
$$\frac{1}{(1-x)(1-x^2)(1-x^3)}=1+x+2x^2+3x^3+4x^4+5x^5+7x^6+8x^7+10x^8\cdots$$
so the rank of $\mathcal{C}(\g(G_2),q,14)$ is $10$.
\end{enumerate}
\end{example}
We close with some remarks.
\begin{remark}
For some Lie types the pre-modular categories described here admit
pre-modular subcategories.   For $\g$ of Lie type $A_r$ and $\ell$
is chosen so that $\gcd(\ell,r+1)=1$, one obtains a modular
subcategory generated by the simple object whose labels are integer
weights (see \cite{MaW}).  The rank of this subcategory is easily
obtained as $\frac{1}{r+1}$ times the rank of the full category.
Similarly ene obtains modular subcategories from $\g$ of type $B_r$
for $\ell$ odd by taking the simple objects labeled by integer
(non-spin) weights give a pre-modular subcategory with rank exactly
half that of the original category.
\end{remark}

\bibliographystyle{amsalpha}

\end{document}